%% file: main_merged.tex
\documentclass[11pt]{article}

\input{macros_Natalia}

\title{A Note on the NP-Hardness on Quadratic Optimization over Doubly Stochastic Matrices and the Gromov-Wasserstein Problem}

\author{Yong Sheng Soh}

\begin{document}

\maketitle

\begin{abstract}
The Gromov-Wasserstein problem is an extension of the regular optimal transport problem that is adapted for pairs of distributions residing in incomparable spaces.  Broadly speaking, given a pair of metric measure spaces, the Gromov-Wasserstein problem seeks couplings between these spaces that minimize the distortion between the metric evaluations.  In this note we record the fact that the Gromov-Wasserstein distance between finite metric spaces is NP-Hard to compute.  Our argument relies on a reduction from the Hamiltonian Cycle Problem, and it builds on ideas in Vayer \cite{vayer2020} and Ryner et. al. \cite{ryner2023} that draw connections between the 2-Gromov-Wasserstein and the Quadratic Assignment Problem.  During the preparation of this manuscript we learned of the appearance of \cite{BGP:26}, whose result essentially supersedes ours.  We hope this note serves as a useful companion whose proof is, hopefully, more elementary and easier to understand.
\end{abstract}

\section{Introduction}

The optimal transport (OT) problem concerns moving mass from one location to another in a manner that minimizes the total cost associated with transportation.  The OT task was originally developed in the context of logistics and economics.  In recent years, the OT task has been applied in a range of scientific and engineering applications including statistics, signal processing, and computer vision tasks among others; today, it forms the backbone of a number of machine learning techniques. 

Formally, the Kantorovich formulation of the optimal transportation problem is expressed as the solution to the following optimization instance:
\begin{equation} \label{eq:kantorovich}
    \min_{\pi} ~~  \int_{X \times Y} c(x,y) d \pi(x,y) \quad
    \mathrm{s.t.} \quad \int_{X} d\pi(x,y) = \mu_Y, \int_{Y} d\pi(x,y) = \mu_X.
\end{equation}
Here, $X$ and $Y$ models the source and the destination space, while $\mu_{X} \in \mathcal{P}(X)$ and $\mu_{Y} \in \mathcal{P}(Y)$ are probability distributions that model the distribution of mass in the source as well as in the destination.  The decision variable $\pi$ is a joint probability distribution over the space $X \times Y$, and it models the amount of mass sent from location $x\in X$ to $y \in Y$.  The {\em coupling constraints} ensure that the transport plan $\pi$ moves mass from the source to the destination completely.

{\bf Optimal Transport between distributions residing in incomparable spaces.} The Kantorovich formulation of the OT problem \eqref{eq:kantorovich} is premised on that one has access to a cost function $c(x,y)$.  This, however, is not possible in certain settings where, for instance, the distributions in the source and destination reside on spaces that are fundamentally different objects, and for which there is no natural notion of cost function.  

To this end, the Gromov-Wasserstein (GW) problem is an alternative formulation that seeks couplings between spaces that minimize the distortion between the metric spaces.  More formally, a {\em metric measure space} is a triplet $(X,d,\mu)$ where $(X,d)$ specifies a metric space, and where $\mu$ specifies a measure over $X$.  Given a pair of metric measure spaces $(X,d_X,\mu_X)$ and $(Y,d_Y,\mu_Y)$, the Gromov-Wasserstein problem seeks a solution to the following problem \cite{Memoli2007, Memoli} :
\begin{equation} \label{eq:gw3}
\begin{aligned}
    \min_{\pi} ~~ & \frac{1}{2} \left( \int_{(X \times Y)^2} \Gamma(x_0,y_0,x_1,y_1)^p d \pi(x_0,y_0) d\pi(x_1,y_1) \right)^{1/p} \\
    \mathrm{s.t.} ~~ & \int_{X} d\pi(x,y) = \mu_Y, \int_{Y} d\pi(x,y) = \mu_X
\end{aligned}.
\end{equation}
Here, the cost $\Gamma$ is typically chosen as the difference between metric evaluations (for some choice of $q$)
\begin{equation*}
    \Gamma(x_0,y_0,x_1,y_1) = | d^q_X(x_0,x_1) - d^q_Y(y_0,y_1) |.
\end{equation*}
The basic intuition behind \eqref{eq:gw3} is that it seeks couplings $\pi$ that minimizes the distortion of the metric evaluations from $X$ to $Y$. 

The GW problem finds applications in a range of computational tasks that require alignment, comparison, or integration of datasets that are not naturally embedded in a common space.  Examples of applications of the GW problem in the biomedical sciences include the Cryo-EM alignment problem \cite{Riahi}, inferring communication between living cells \cite{Cang}, integrating data from different experiments \cite{Demetci}, classifying biomedical time series \cite{Kravtsova}, as well as quantifying evolution of gene expression networks \cite{oliver2026}.  Examples of applications of the GW problem in machine learning include generative modeling tasks \cite{Bunne}, reinforcement learning \cite{Fickinger}, graph partitioning and clustering \cite{ChowdhuryNeedham, Vayer, Xu2022}, and graph classification tasks \cite{hohmeier2026}.  Examples of applications of the GW problem in the statistical sciences include developing methodology for testing similarity \cite{ steinkamp2026,weitkamp2024}, quantifying convergence of estimators \cite{han2025,otcrossdomain2025,wu2026,zhang2024}, and feature selection tasks \cite{LYLC:26}.

\subsection{Computing the Gromov-Wasserstein distance}  

In settings where the source and target spaces are over a finite number of points, the GW problem \eqref{eq:gw3} can be cast as the following finite dimensional optimization instance
\begin{equation} \label{eq:gw_discrete}
\begin{aligned}
    \min_{\pi} ~~ & \frac{1}{2} \Big( \sum_{ijkl} \Gamma(x_i,x_j,y_k,y_l)^{p}  \pi_{ik} \pi_{jl} \Big)^{1/p} \\
    \mathrm{s.t.} ~~ & \pi \mathbf{1} = \mu_X, \pi^T \mathbf{1} = \mu_Y, \pi_{ij} \geq 0 ~ \forall i,j
\end{aligned}.
\end{equation}
The decision variable $\pi$ is now modeled as a finite dimensional matrix.  Here, we use $\mathbf{1}$ to denote the vector whose entries are all ones.  We make two further simplifications.  First, we raise the objective to the $p$-th power to eliminate the appearance of the $p$-th root.  Second, we note that the metric information $d_X$ and $d_Y$ can be conveyed via finite dimensional matrices $P \in \mathbb{R}^{|\mu_X| \times |\mu_X|}$ and $Q \in \mathbb{R}^{|\mu_Y| \times |\mu_Y|}$ whereby
\begin{equation} \label{eq:q-here}
P_{ij} = d^q_{X} (x_i,x_j) \quad \text{and} \quad Q_{kl} = d^q_{Y} (y_k,y_l).    
\end{equation}
By applying these simplifications, the formulation \eqref{eq:gw_discrete} is equivalent to the following quadratic programming instance
\begin{equation} \label{eq:gw_discrete-2}
\begin{aligned}
    \min_{\pi} ~~ & \sum_{ijkl} | P_{ij} - Q_{kl} |^{p}  \pi_{ik} \pi_{jl}  \\
    \mathrm{s.t.} ~~ & \pi \mathbf{1} = \mu_X, \pi^T \mathbf{1} = \mu_Y, \pi_{ij} \geq 0 ~ \forall i,j.
\end{aligned}
\end{equation}
In particular, the objective is quadratic in $\pi$ while the constraints are affine.

{\bf Computational complexity of \eqref{eq:gw_discrete-2}.}  The central question of this note is to address the computational complexity of \eqref{eq:gw_discrete-2}.  The general class of quadratic programs is NP-hard \cite{Vavasis2009}.  Unfortunately, the choice of $\Gamma$ specified as the difference of distance matrices as in \eqref{eq:gw_discrete-2} leads to a non-convex objective for any input data \cite{kravtsova2026note}.  Vayer, in his thesis \cite{vayer2020}, as well as Ryner et. al. in subsequent work \cite{ryner2023}, show that the GW instantiated with $p=2$ and with cost matrices $P$ and $Q$ arising from Euclidean point clouds is equivalent to a Quadratic Assignment Problem (QAP) instance.  To this end, the general class of QAPs is NP-Hard.  These connections suggest that \eqref{eq:gw_discrete-2} is computational intractable, and is a point raised by many prior works -- see, for example, the discussion in Section 7 of \cite{Memoli}, Section 1.1 of \cite{Peyre}, Section 1 of \cite{Scetbon}, and/or Section 3 of \cite{Chen}.  

The precise version of the type of result we seek in this paper is a formal reduction of a known class of NP-Hard problems to an optimization formulation of the form \eqref{eq:gw_discrete} or \eqref{eq:gw_discrete-2}, and where $P$ and $Q$ are distance matrices that arise from a {\em bona fide} metric space.  We were initially not aware of a existing result to this effect.  The goal of writing this note is to record such a fact.  However, as we were about to post this note, it came to our attention that the Breitung et. al. independently showed NP-Hardness of the GW problem \cite{BGP:26}.  The result in \cite{BGP:26} is more general than ours -- first, it discusses NP-Hardness in a broader class of problems concerning distances over graphs, and within the context of the GW problem it applies to general $p$ whereas our result only applies to $p=2$.  In our opinion, the proof supplied in this work is more elementary, and quite possibly easier to understand as it (essentially) disguises a Hamiltonian Cycle instance into a GW instance.  We certainly hope the ideas in our note will prove to be valuable alongside the work in \cite{BGP:26}.

{\bf Practical algorithms for solving \eqref{eq:gw_discrete-2}.}  We briefly discuss existing algorithms for solving the GW problem.  These include a Conditional Gradient / Frank-Wolfe based method \cite{vincent-cuaz21a}, one based on entropic regularization \cite{Peyre}, one based on low-rank assumptions on the cost \cite{Scetbon} (among many others).  These methods are currently deployed on Python toolboxes \cite{CuturiOTbox,POT}.  Unfortunately, these methods are only able to compute local or approximate solutions, and are generally unable to provide convergence guarantees to a global solution.  There is a different line of work based on semidefinite relaxations -- these have the advantage in that they do solve \eqref{eq:gw_discrete-2} to global optimality frequently \cite{Chen,TNS:26,TS:26} while providing a certificate of global optimality; however, these methods require the solution of a semidefinite programming instance, which is considerably more expensive compared to local methods.  On this note, a closely related approach is one based on the moment--SOS hierarchy \cite{MN:24,TS:26} -- these are capable of computing the optimal value to the GW problem in a fashion that is independent of how the metric space is discretized; the drawback is that it is often challenging to extract the corresponding transport map.  Last, we mention analytical approaches -- these include computing exact values for specific families of distributions \cite{AAELMP:26,SDD:21} and lower bounds as in \cite{Memoli}, as well as adaptations for networks \cite{Chowdhury} and ultrametric spaces \cite{Memoliultra}.  These techniques were applied in \cite{kravtsova2025,Le,steinkamp2026,weitkamp2024}.

\subsection{Quadratic Optimization over Doubly Stochastic Matrices and Related Problems}

Our reduction will focus on the specific setting where (i) the metric spaces are finite collection of an equal number of points which we denote by $n$, and (ii) the associated measures are taken to be uniform.  In this case, by writing $X := n\pi$, the optimization instance \eqref{eq:gw_discrete-2} is equivalent to the following form:
\begin{equation} \label{eq:quad_over_ds-2}
\begin{aligned}
    \underset{X \in \mathbb{R}^{n \times n}}{\min} \quad & \sum_{ijkl} | P_{ij} - Q_{kl} |^{p} X_{ik} X_{jl} \\
    \mathrm{s.t.} \quad & X \mathbf{1} = \mathbf{1}, X^T \mathbf{1} = \mathbf{1}, X_{ij} \geq 0 ~ \forall i,j
\end{aligned}.
\end{equation}
In particular, the feasible region is the collection of doubly stochastic matrices:
$$
\ds(n) := \{ X : X \mathbf{1} = \mathbf{1}, X^T \mathbf{1} = \mathbf{1}, X \geq 0 \} \subset \mathbb{R}^{n \times n}.
$$

The form of the problem \eqref{eq:quad_over_ds} bears resemblance to a number of other optimization instances that are known to be NP-Hard, and so we give a quick discussion of these problems.  In what follows we let $\Delta(n)$ denote the probability simplex over $\mathbb{R}^{n}$:
$$
\Delta(n) := \{ \bx : \bx^T \mathbf{1} = 1, \bx \geq 0 \} \subset \mathbb{R}^{n}.
$$
We let $\perm(n)$ denote the collection of permutation matrices over $\mathbb{R}^{n \times n}$:
$$
\perm(n) := \{ X : X \mathbf{1} = \mathbf{1}, X^T \mathbf{1} = \mathbf{1}, X_{ij} \in \{0,1\} \, \forall (i,j) \} \subset \mathbb{R}^{n \times n}.
$$
Recall that from the Birkhoff-von Neumann's theorem that $\ds(n)$ is the convex hull of $\perm(n)$; moreover, the extreme points of $\ds(n)$ are $\perm(n)$ itself.  Last, we let $I$ denote the identity matrix, and $J$ be the matrix whose entries are all ones.

{\bf Quadratic Optimization over the Probability Simplex.}  Consider the problem of minimizing a quadratic objective over the probability simplex
\begin{equation} \label{eq:quad_over_simplex}
    \underset{\bx \in \mathbb{R}^{n}}{\min} \quad \bx^T C \bx \quad
    \mathrm{s.t.} \quad \bx \in \Delta(n).
\end{equation}
Given a doubly stochastic matrix $X\in \ds(n)$, the sum of its entries equals to $n$, which is constant.  Hence the set of doubly stochastic matrices is a strict subset of the probability simplex.  In this sense, \eqref{eq:quad_over_ds} can be viewed as an optimization instance of the type \eqref{eq:quad_over_simplex}, but with additional constraints.

Let $\mathcal{G} = (\mathcal{V},\mathcal{E})$ be a graph.  Let $n=|\mathcal{V}|$ be the number of vertices in $\mathcal{G}$.  Let $\alpha(\mathcal{G})$ be the {\em stability number} of the graph $\mathcal{G}$.  It is known from Motzkin-Strauss \cite{MS:65} that
\begin{equation*}
    \frac{1}{\alpha(\mathcal{G})} ~~ = ~~ \min_{\bx} \quad \bx^T (A_{\mathcal{G}}+I_n) \bx \quad
    \mathrm{s.t.} \quad \bx \in \Delta(n).
\end{equation*}
Here, $A_{\mathcal{G}}$ is the adjacency graph of $\mathcal{G}$.  The stability number of a graph is NP-Hard to compute, and hence the class of problems described in \eqref{eq:quad_over_simplex} is also NP-Hard.

{\bf Quadratic Assignment Problem.}  The Quadratic Assignment Problem (QAP) is given as follows:  Given matrices $A$ and $B$, find the permutation matrix $X$ minimizing the following instance:
\begin{equation} \label{eq:qap}
\begin{aligned}
    \underset{X \in \mathbb{R}^{n \times n}}{\min} \quad & \mathrm{tr} (A X B X^T) \quad \mathrm{s.t.} \quad X \in \perm(n).
\end{aligned}.
\end{equation}
The class of QAPs is known to be NP-Hard.  For instance, every Hamiltonian Cycle instance can be encoded as a QAP instance.  As we discuss later, in the setting where the cost $\Gamma$ in \eqref{eq:gw3} is chosen to be a squared-loss, the objective reduces to a QAP-type loss as in \eqref{eq:qap}; the difference is that \eqref{eq:gw3} takes place over $\ds(n)$ whereas \eqref{eq:qap} takes place over $\perm(n)$.

{\bf The Hamiltonian Cycle Problem.}  Let $\mathcal{G} = (\mathcal{V},\mathcal{E})$ be an undirected graph.  A Hamiltonian cycle is a cycle where consecutive vertices $v_i v_{i+1} \in \mathcal{E}$ are edges, and where each vertex is visited exactly once.  Given a graph $\mathcal{G}$, the Hamiltonian Cycle Problem decides if there exists a Hamiltonian Cycle.  The Hamiltonian Cycle Problem is NP-Hard, and is one of Karp's list of 21 problems.  There is a simple reduction of the Hamiltonian Cycle instance to a QAP instance: Let the matrices $A$ and $B$ each encode a cycle and a generic graph.  Then solving the associated QAP instance is equivalent to finding a Hamiltonian Cycle within the graph encoded in $B$.  We rely on a variant of this reduction later in our analysis.

\section{Quadratic Optimization over Doubly Stochastic Matrices}

Consider the following class of quadratic optimization over the set of doubly stochastic matrices:
\begin{equation} \label{eq:quad_over_ds}
\begin{aligned}
    \underset{X \in \mathbb{R}^{n \times n}}{\min} \quad & \sum_{ijkl} C_{ijkl} X_{ik} X_{jl} \\
    \mathrm{s.t.} \quad & X \mathbf{1} = \mathbf{1}, X^T \mathbf{1} = \mathbf{1}, X_{ij} \geq 0 ~ \forall i,j
\end{aligned}.
\end{equation}
The difference between \eqref{eq:quad_over_ds} and \eqref{eq:quad_over_ds-2} is that the cost $C_{ijkl}$ need not be a difference of distance matrices.  

\begin{theorem} \label{thm:reduction_to_qap}
    There is a polynomial-time reduction of QAPs to the problem class \eqref{eq:quad_over_ds} (of the same problem size).  
\end{theorem}

As a corollary, the problem class \eqref{eq:quad_over_ds} is NP-Hard.  Note this conclusion will be superseded by our discussion in the next section.  The purpose of stating Theorem \ref{thm:reduction_to_qap} is that it has an elementary proof.  The proof relies on the fact that the maximum of a convex function over a compact convex set is attained at its extreme points, a result sometimes referred to as Bauer's maximum principle:
\begin{proposition}[Bauer's Maximum Principle] \label{thm:bauer}
    Let $C \subset \mathbb{R}^{n}$ be compact convex.  Let $f$ be convex and continuous.  The $f$ attains its maximum over $C$ at an extreme point of $C$.
\end{proposition}

\begin{proof} [Proof of Theorem \ref{thm:reduction_to_qap}]
    Begin with a QAP instance specified by the matrices $A,B \in \mathbb{R}^{n \times n}$.  Note that $\mathrm{tr}(AXBX^T) = \sum_{ijkl} A_{ji} B_{kl} X_{ik} X_{jl}$.  Define 
    $$
    C_{ijkl} = A_{ji} B_{kl} - \begin{cases}
        \gamma \quad \text{ if } (i,j) = (k,l) \\
        0 \quad \text{ otherwise }    
    \end{cases},
    $$
    where $\gamma$ is a constant to be specified later.  The objective in \eqref{eq:quad_over_ds} can then be re-written as
    $$
    f(X) ~=~ \sum_{ijkl} C_{ijkl} X_{ik} X_{jl} 
    =~ \sum_{ijkl} A_{ji} B_{kl} X_{ik} X_{jl} - \gamma \| X \|_F^2.
    $$
    Note here that $\|\cdot \|_F^2$ is the Euclidean-norm over $\mathbb{R}^{n \times n}$.  We pick $\gamma$ sufficiently large so that $f$ is strictly concave.  

    Return to the optimization instance \eqref{eq:quad_over_ds}.  Because we are minimizing a concave function over a polyhedral constraint set, the optimal solution is attained at an extreme point (Proposition \ref{thm:bauer}).  By Birkhoff-von Neumann's theorem, the extreme points of $\ds(n)$ are $\perm(n)$.  Hence the optimization instance
    $$
    \min_{X} ~~ f(X) \quad \mathrm{s.t.} \quad X \in \ds(n)
    $$
    is equivalent to 
    $$
    \min_{X} ~~ f(X) \quad \mathrm{s.t.} \quad X \in \perm(n)
    $$
    Now note that $\| X \|_F^2 = n$ for all $X \in \perm(n)$.  Hence, by subtracting a constant $n$ from the objective, the above is equivalent to
    $$
    \min_{X} ~~ \sum_{ijkl} A_{ji} B_{kl} X_{ij} X_{kl} \quad \mathrm{s.t.} \quad X \in \perm(n).
    $$
    This is exactly the QAP formulation.

    Last, we note that problem transformation operates in the same dimensions, and hence these associated reductions are polynomial.
\end{proof}

\section{The Gromov-Wasserstein Problem}

In this section we shift our focus to the Gromov-Wasserstein problem.  The goal of this section is to show a reduction of a Hamiltonian Cycle instance to a Gromov-Wasserstein instance  \eqref{eq:gw3}.  

As a reminder, the Hamiltonian Cycle problem is the following decision problem:  
\begin{quote}
    Given a graph $\mathcal{G} = (\mathcal{V},\mathcal{E})$, decide if there exists a Hamiltonian Cycle that visits every vertex exactly once, and returns to the origin.    
\end{quote}

{\bf [Step 1 -- Specifying the matrices $P$ and $Q$ and the corresponding mm-spaces]:}  Let $n = |\mathcal{V}|$.  Let $\mathcal{C}:=\mathcal{C}_n$ be the cyclic graph over $n$ vertices.  Define the matrices 
$$
P := 2(J-I) - A_{\mathcal{C}}, \quad Q := 2(J-I) - \epsilon A_{\mathcal{G}}.
$$
Here, $A_{\mathcal{G}}$ is the adjacency graph associated with $\mathcal{G}$.  Also, $\epsilon$ is chosen such that $0 < \epsilon \leq 1/\|A_{\mathcal{G}}\|$, where $\|\cdot\|$ denotes the spectral norm of the input matrix.

The associated metric measure spaces $(X,d_X,\mu_X)$ and $(Y,d_Y,\mu_Y)$ where (i) $X$ and $Y$ are a discrete collection of $n$ points, (ii) $d_X(x_i,x_j) = \sqrt{P_{ij}}$ and $d_Y(y_k,y_l) = \sqrt{Q_{kl}}$, and (iii) $\mu_X$ and $\mu_Y$ are the uniform measures.  (In particular, this means that we pick $q=2$ inside \eqref{eq:q-here}.)

We briefly explain why the choices $d_X$ and $d_Y$ specify an actual metric space.  We say that a matrix $D \in \mathbb{R}^{n\times n}$ is an {\em Euclidean Distance Matrix} if $D_{ij} = \| \bx_i - \bx_j \|_2^2$ for some collection of vectors $\{\bx_i\}_{i=1}^{n}$ in some Euclidean space.  Suppose we define
$$
H:= I - J/n, 
$$
and we define
$$
\hat{D} = HDH.
$$
Then $D$ is an Euclidean Distance Matrix if and only if $\hat{D}$ is negative semidefinite -- this is sometimes referred to as the Schoenberg criterion.  In particular:

\begin{proposition} \label{thm:pq-negdef}
    We have $\hat{P} := HPH \preceq 0$ and $\hat{Q} := HQH \preceq 0$.
\end{proposition}

\begin{proof}[Proof of Proposition \ref{thm:pq-negdef}]
First, note that $H$ is a projection matrix that projects out components in $\boldsymbol{1}/\sqrt{n}$.  In particular, because $H$ is a projection matrix, one has $H^2 = H$.  Second, a direct calculation shows that
$$
-H ( 2(J-I) - A) H = H (2I + A) H.
$$
    
[$-\hat{P}\succeq 0$]:  Using standard facts about the adjacency graph of a cycle, the eigenvalues of the adjacency matrix $A_{\mathcal{C}}$ is $\{ 2 \cos (2\pi k /n) : 0 \leq k \leq n-1 \}$.  Subsequently, the eigenvalues of $2I+A_{\mathcal{C}}$ are $\{ 2+2 \cos (2\pi k /n) : 0 \leq k \leq n-1 \}$.  These are non-negative, which implies $2I + A_{\mathcal{C}} \succeq 0$ and subsequently $-\hat{P} \succeq 0$.
    
[$-\hat{Q} \succeq 0$]: Earlier we noted that $- \hat{Q} = H (2I + \epsilon A_{\mathcal{G}} ) H$.  Since $\epsilon \leq 2/\|A_{\mathcal{G}}\|$, it follows that all the eigenvalues of $2I + \epsilon A_{\mathcal{G}}$ are non-negative.  Subsequently $2I + \epsilon A_{\mathcal{G}} \succeq 0$, which then implies $- \hat{Q} \succeq 0$.
\end{proof}

{\bf [Step 2 -- Show that the objective coincides with a convex function over $\ds(n)$]:}  Our next step is to show a formal reduction to a specific instance of the GW instance as in \eqref{eq:gw_discrete-2}.  In the above, our choice of mm-spaces is such that the source and target have an equal number of points, and where the source and target measures are uniform.  As such, the optimization instance \eqref{eq:gw_discrete-2} simplifies to \eqref{eq:quad_over_ds-2}.  In the remainder of this note, we pick $p=2$.  Then, by expanding the squared loss objective, we obtain
$$
(P_{ij} - Q_{kl})^2 = P_{ij}^2 + Q_{kl}^2 - P_{ij} Q_{kl}.
$$
For all $X \in \ds(n)$ one has
$$
\sum_{ijkl} P_{ij}^2 X_{ik} X_{jl} = \Big(\sum_{ij} P_{ij}^2 \big(\sum_{k} X_{ik} \big) \big(\sum_l X_{jl} \big) \Big) = \sum_{ij} P_{ij}^2 = \| P \|_F^2, 
$$
which is a constant independent of $X$.  In particular, by eliminating the constants, the optimization instance \eqref{eq:quad_over_ds-2} simplifies to
\begin{equation} \label{eq:gw3}
    \max_{X \in \mathbb{R}^{n \times n} } ~~ \sum_{ijkl} P_{ij}Q_{kl} X_{ik} X_{jl} \quad
    \mathrm{s.t.} \quad X \in \ds(n).
\end{equation}

The objective is given by
$$
g(X) := \sum_{ijkl} P_{ij}Q_{kl} X_{ik} X_{jl} = \mathrm{tr}(PXQX^T).
$$
Consider the function 
$$
    \hat{g}(X) := \mathrm{tr}(\hat{P} X \hat{Q} X^T).
$$

\begin{lemma} \label{thm:is_ghat}
    The function $\hat{g}(X)$ coincides with $g(X) + l(X) + c$ over $X \in \ds(n)$.  Here,  $l(X)$ is a linear function and $c$ is a constant.
\end{lemma}

\begin{proof}
First, note that for any $X \in \ds(n)$ one has $XJ = JX = J$.  This gives $HXH = X - J/n$ and $HX^TH = X^T - J/n$.  Hence, for any $X \in \ds(n)$,
\begin{equation*}
\begin{aligned}
    & \mathrm{tr}\bigl( (HPH)X(HQH)X^T \bigr)  \\
    = ~ & \mathrm{tr}\bigl( P(HXH) Q (HX^TH) \bigr) \\
    = ~ & \mathrm{tr}\bigl( P(X - J/n) Q (X^T - J/n) \bigr) \\
    = ~ & \mathrm{tr}(PXQX^T) - \frac{1}{n}\mathrm{tr}(PXQJ) - \frac{1}{n}\mathrm{tr}(PJQX^T) + \frac{1}{n^2} \mathrm{tr}(PJQJ)
\end{aligned}.
\end{equation*}
In particular, $l(X) = - \mathrm{tr}(PXQJ)/n - \mathrm{tr}(PJQX^T)/n$ is linear in $X$ while $c = \mathrm{tr}(PJQJ)/n^2$ is a constant.
\end{proof}

As a reminder, the function $\hat{g}$ can be written in terms of Kronecker products.  Let $\mathrm{vec}(X)$ denote the vectorization of the matrix $X$.  Then
$$
    \hat{g}(X) = \mathrm{tr}(\hat{P} X \hat{Q} X^T) = \langle \mathrm{vec}(X), \mathrm{vec}(\hat{P} X \hat{Q}) \rangle = \langle \mathrm{vec}(X), [\hat{Q} \otimes \hat{P}] \mathrm{vec}(X) \rangle.
$$
Note that in the last line we used the fact that $\hat{Q}$ is symmetric and hence $\hat{Q}^T = \hat{Q}$.

\begin{lemma} \label{thm:g_convex} 
    The function $\hat{g}$ is convex.
\end{lemma}

\begin{proof}  Recall a property regarding Kronecker products that the eigenvalues of $\hat{Q} \otimes \hat{P}$ are the product of the eigenvalues of $\hat{Q}$ and $\hat{P}$.  Subsequently the result follows from Proposition \ref{thm:pq-negdef}.
\end{proof}

These results allow us to arrive at the following conclusion.  Consider the optimization instances 
\begin{equation} \label{eq:red-1}
    \max_{X} ~~ \mathrm{tr} (P X Q X^T) \quad \mathrm{s.t.} \quad X \in \ds(n),    
\end{equation}
as well as 
\begin{equation} \label{eq:red-2}
    \max_{X} ~~ \mathrm{tr} (P X Q X^T) \quad \mathrm{s.t.} \quad X \in \perm(n).    
\end{equation}

\begin{proposition}
    The optimization instances \eqref{eq:red-1} and \eqref{eq:red-2} are equivalent.
\end{proposition}

\begin{proof}
    By Lemma \ref{thm:is_ghat}, \eqref{eq:red-1} is equivalent to 
\begin{equation} \label{eq:red-3}
    \max_{X} ~~ \mathrm{tr} (\hat{P} X \hat{Q} X^T) \quad \mathrm{s.t.} \quad X \in \ds(n).    
\end{equation}
From Lemma \ref{thm:g_convex} the objective is convex.  Hence the optimal solution is attained at the extreme points of $\ds(n)$ (Proposition \ref{thm:bauer}).  By Birhkoff-von Neumann's theorem the set of extreme points of $\ds(n)$ is $\perm(n)$.  Therefore \eqref{eq:red-3} is equivalent to 
\begin{equation*}
    \max_{X} ~~ \mathrm{tr} (\hat{P} X \hat{Q} X^T) \quad \mathrm{s.t.} \quad X \in \perm(n).
\end{equation*}
By Lemma \ref{thm:is_ghat} once again this is equivalent to \eqref{eq:red-2}.
\end{proof}

{\bf [Step 3 -- Reducing Hamiltonian Cycle to GW]:}  
Given $X \in \perm(n)$, the objective simplifies to
$$
\mathrm{tr}(PXQX^T) = \sum_{i,j} P_{ij} Q_{\sigma(i)\sigma(j)}.
$$
Using the definition $P = 2(J-I) - A_{\mathcal{C}}$, the objective further simplifies to 
\begin{equation} \label{eq:Q_sigma}
    \sum_{i,j} P_{ij} Q_{\sigma(i)\sigma(j)} = 2 \Bigl[\sum_{k \neq l} Q_{kl}- \Bigl( \sum_{i=1}^{n-1} Q_{\sigma(i)\sigma(i+1)} + Q_{\sigma(1)\sigma(n)}\Bigr) \Bigr].
\end{equation}
The term $\sum_{k \neq l} Q_{kl}$ is independent of $\sigma$.  Therefore the optimization instance \eqref{eq:gw3} is equivalent to
\begin{equation} \label{eq:red-4}
    \min_{\sigma:[n] \rightarrow [n]} \sum_{i=1}^{n-1} Q_{\sigma(i)\sigma(i+1)} + Q_{\sigma(1)\sigma(n)}.
\end{equation}
Recall that our choice of $Q = 2(J-I) - A_{\mathcal{G}}$ gives 
\begin{equation*}
    Q_{ij} = \begin{cases}
    2- \frac{1}{n},  & \text{ if vertex } i \text{ is connected to vertex } j \text{ in } \mathcal{G} \\
    2, & \text{ else}
\end{cases} 
\end{equation*}
Therefore the optimal value to \eqref{eq:red-4} is
\begin{equation*}
\text{Opt Value of } \eqref{eq:red-4} \quad = \quad
 \begin{cases}
  n(2 - \frac{1}{n}) & \text{ if } \mathcal{G} \text{ has a Hamiltonian cycle } \\
 \geq  n(2 - \frac{1}{n}) +1/n &  \text{ if }\mathcal{G} \text{ does not have Hamiltonian cycle}
\end{cases} 
\end{equation*}

In particular if we denote
\begin{equation}\label{eq:t}
t(\mathcal{G}):=\|P\|_F^2 + \|Q\|_F^2 -  4 \Bigl[\sum_{k\neq l} Q_{kl} - n(2 - \frac{1}{n}) \Bigr],
\end{equation}
then
\begin{equation}\label{eq:cyclic_sum}
\mathrm{GW}_2^2((X,d_X,\mu_X), (Y,d_Y,\mu_Y)) \quad = \quad
 \begin{cases}
 = t(\mathcal{G}) & \text{if $\mathcal{G}$ has a Hamiltonian cycle} \\
 
 \geq t(\mathcal{G}) + 1/n & \text{if $\mathcal{G}$ does not have Hamiltonian cycle}
\end{cases} 
\end{equation}

 
Last we note that our transformation from Hamiltonian Cycle Problem to \eqref{eq:gw3} can be executed in polynomial time in terms of the size of the graph $\mathcal{G}$.  Indeed, to perform our transformation, we need to construct distance matrices $P$ and $Q$ by placing values $\{0,1,2\}$ in $P$ and $\{0,2-\frac{1}{n},2\}$ in $Q$ according to whether a corresponding pair of vertices from the cycle graph $\mathcal{C} _n$ and the graph of interest $\mathcal{G} $ has an edge -- this can be done in $O(n^2)$.  In addition, computing the threshold value $t(\mathcal{G})$ can be done in $O(n^2)$.  In particular:

\begin{theorem} \label{thm:GW_NPhard}
    There is a polynomial time reduction from Hamiltonian Cycle to the family of optimization instances specified by \eqref{eq:gw3} where $P$ and $Q$ are distance matrices. 
\end{theorem}

\subsection{Relation to prior work}

The reduction presented here depends on a number of observations that have appeared in prior work.  In particular,  \cite[Proposition 5.5.1]{vayer2020} as well as \cite[Proposition 1 and 2]{ryner2023} explain that the GW problem, specialized to $p=q=2$ and distance matrices arising from Euclidean point clouds, give rise to a quadratic program of the form \eqref{eq:gw3} where the objective is convex.  Furthermore, \cite[Proposition 2]{ryner2023} explain that the matrices $P$ and $Q$ may be chosen suitably so that the resulting objective is convex, in which case we may reduce \eqref{eq:gw3} to one over the set of permutation matrices.  As such, the real intellectual contribution is that we show that it is possible to encode the Hamiltonian Cycle problem using choices of $P$ and $Q$.

\section{Acknowledgments}

The author wishes to acknowledge discussion with Natalia Kravtsova.  In fact, an earlier version of this note was written in collaboration with Natalia.  The author wishes to acknowledge helpful discussions with Facundo Mémoli and Shmuel Friedland.


\section{AI Disclosure}

As a matter of full disclosure, the reduction of the Hamiltonian Cycle Problem to \eqref{eq:gw3} was discovered with assistance from the Claude LLM, using the Opus 5 model running on Max setting.  We prompted the LLM to discover a reduction from the Traveling Salesman Problem (TSP) to the GW problem formulation as it is stated in \eqref{eq:gw3}.  The response is a reduction from Hamiltonian Cycle Problem, which is (essentially) similar to the TSP.  

\bibliographystyle{plain}
\bibliography{gw_nphard}

\appendix

\end{document}

%% file: macros_Natalia.tex
\usepackage{amssymb,amsthm,amsmath}
\usepackage[pdftex]{graphicx}
\usepackage[margin=1in]{geometry}
\usepackage{bm}

\usepackage[colorlinks=true,linkcolor=blue,citecolor=blue]{hyperref}

\newtheorem{theorem}{Theorem}[section]
\newtheorem{lemma}[theorem]{Lemma}
\newtheorem{proposition}[theorem]{Proposition}

\theoremstyle{definition}

\theoremstyle{remark}

\newcommand{\bx}{\mathbf{x}}

\newcommand{\perm}{\mathrm{Perm}}
\newcommand{\ds}{\mathrm{DS}}

%% file: gw_nphard.bib
@article{BGP:26,
  title={Graph Distortion Distances, Gromov–Hausdorff and the Hardness of Gromov–Wasserstein},
  author={Jens Breitung and Timo Gervens and Julius P\"{o}rzgen},
  journal={arXiv preprint  	arXiv:2610.1692},
  year={2026}
}

@article{AAELMP:26,
  title={The Gromov–Wasserstein Distance Between Spheres},
  author={Shreya Aryaand and Arnab Auddy and Ranthony Edmonds and Sunhyuk Lim and Facundo Memoli and Daniel Packer},
  journal={Foundations of Computational Mathematics},
  year={2026}
}

@article{SDD:21,
  title={Gromov-Wasserstein Distances between Gaussian Distributions},
  author={Antoine Salmona and Julie Delon and Agnès Desolneux},
  journal={arXiv preprint  	arXiv:2104.07970},
  year={2021}
}

@article{TS:26,
  title={Quadratic Optimization over Probability Measures with Coupling Constraints},
  author={Hoang Anh Tran and Yong Sheng Soh},
  journal={arXiv preprint  arXiv:2608.24184},
  year={2026}
}

@article{TNS:26,
  title={Sum-of-{S}quares {H}ierarchy for the {G}romov--{W}asserstein {P}roblem},
  author={Hoang Anh Tran and Binh Tuan Nguyen and Yong Sheng Soh},
  journal={SIAM Journal on Optimization},
  volume={36},
  number={2},
  pages={729--759},
  year={2026},
  publisher={SIAM}
}

@article{MN:24,
 Author = {Mula, Olga and Nouy, Anthony},
 Title = {Moment-{SoS} {M}ethods for {O}ptimal {T}ransport {P}roblems},
 FJournal = {Numerische Mathematik},
 Journal = {Numer. Math.},
 ISSN = {0029-599X},
 Volume = {156},
 Number = {4},
 Pages = {1541--1578},
 Year = {2024},
 Language = {English},
 DOI = {10.1007/s00211-024-01422-x},
 zbMATH = {7896230},
 Zbl = {1545.49046}
}

@InProceedings{vincent-cuaz21a,
  title = 	 {Online Graph Dictionary Learning},
  author =       {C{\'e}dric Vincent-Cuaz and Titouan Vayer and R{\'e}mi Flamary and Marco Corneli and Nicolas Courty},
  booktitle = 	 {Proceedings of the 38th International Conference on Machine Learning},
  pages = 	 {10564--10574},
  year = 	 {2021},
  editor = 	 {Meila, Marina and Zhang, Tong},
  volume = 	 {139},
  series = 	 {Proceedings of Machine Learning Research},
  month = 	 {18--24 Jul},
  publisher =    {PMLR}
}

@article{LYLC:26,
  title={Fused {G}romov-{W}asserstein {D}istance with {F}eature {S}election},
  author={Harlin Lee and Ying Yu and Mingxin Li and Ranthony Clark},
  journal={arXiv preprint  arXiv:2605.12161},
  year={2026}
}

@article{MS:65,
  title={Maxima for {G}raphs and a {N}ew {P}roof of a {T}heorem of {T}urán},
  author={T. S. Motzkin and E. G. Straus},
  journal={Canadian Journal of Mathematics},
  volume={17},
  year={1965}
}

@phdthesis{V,
  AUTHOR =	 "M. Vajiac",
  TITLE =	 "Gauge theory techniques in quantum cohomology",
  school =	 {Boston University},
  YEAR =	 2000
}

@article{Memoli,
  title={Gromov--{W}asserstein distances and the metric approach to object matching},
  author={M{\'e}moli, Facundo},
  journal={Foundations of computational mathematics},
  volume={11},
  pages={417--487},
  year={2011},
  publisher={Springer}
}

@inproceedings{Memoli2007,
booktitle = {Eurographics Symposium on Point-Based Graphics},
editor = {M. Botsch and R. Pajarola and B. Chen and M. Zwicker},
title = {{On the use of {G}romov--{H}ausdorff distances for shape comparison}},
author = {M{\'e}moli, Facundo},
year = {2007},
publisher = {The Eurographics Association},
ISSN = {1811-7813},
ISBN = {978-3-905673-51-7},
DOI = {/10.2312/SPBG/SPBG07/081-090}
}

@article{POT,
  author  = {R{\'e}mi Flamary and Nicolas Courty and Alexandre Gramfort and Mokhtar Z. Alaya and Aur{\'e}lie Boisbunon and Stanislas Chambon and Laetitia Chapel and Adrien Corenflos and Kilian Fatras and Nemo Fournier and L{\'e}o Gautheron and Nathalie T.H. Gayraud and Hicham Janati and Alain Rakotomamonjy and Ievgen Redko and Antoine Rolet and Antony Schutz and Vivien Seguy and Danica J. Sutherland and Romain Tavenard and Alexander Tong and Titouan Vayer},
  title   = {{POT}: {P}ython {O}ptimal {T}ransport},
  journal = {Journal of Machine Learning Research},
  year    = {2021},
  volume  = {22},
  number  = {78},
  pages   = {1-8},
  url     = {http://jmlr.org/papers/v22/20-451.html}
}

@article{CuturiOTbox,
  title={{O}ptimal {T}ransport {T}ools ({OTT}): A {JAX} toolbox for all things {W}asserstein},
  author={Cuturi, Marco and Meng-Papaxanthos, Laetitia and Tian, Yingtao and Bunne, Charlotte and Davis, Geoff and Teboul, Olivier},
  journal={arXiv preprint arXiv:2201.12324},
  year={2022}
}

@inproceedings{Le,
  title={Flow-based alignment approaches for probability measures in different spaces},
  author={Le, Tam and Ho, Nhat and Yamada, Makoto},
  booktitle={International conference on artificial intelligence and statistics},
  pages={3934--3942},
  year={2021},
  organization={PMLR}
}

@article{Cang,
  title={Inferring spatial and signaling relationships between cells from single cell transcriptomic data},
  author={Cang, Zixuan and Nie, Qing},
  journal={Nature communications},
  volume={11},
  number={1},
  pages={2084},
  year={2020},
  publisher={Nature Publishing Group UK London}
}

@inproceedings{Fickinger,
  title={Cross-Domain Imitation Learning via Optimal Transport},
  author={Fickinger, Arnaud and Cohen, Samuel and Russell, Stuart and Amos, Brandon},
  booktitle={10th International Conference on Learning Representations, ICLR},
  year={2022}
}

@inproceedings{ChowdhuryNeedham,
  title={Generalized spectral clustering via {G}romov--{W}asserstein learning},
  author={Chowdhury, Samir and Needham, Tom},
  booktitle={International Conference on Artificial Intelligence and Statistics},
  pages={712--720},
  year={2021},
  organization={PMLR}
}

@article{Vayer,
  title={Fused {G}romov--{W}asserstein distance for structured objects},
  author={Vayer, Titouan and Chapel, Laetitia and Flamary, R{\'e}mi and Tavenard, Romain and Courty, Nicolas},
  journal={Algorithms},
  volume={13},
  number={9},
  pages={212},
  year={2020},
  publisher={MDPI}
}

@article{Xu2022,
  title={Representing graphs via {G}romov--{W}asserstein factorization},
  author={Xu, Hongteng and Liu, Jiachang and Luo, Dixin and Carin, Lawrence},
  journal={IEEE Transactions on Pattern Analysis and Machine Intelligence},
  volume={45},
  number={1},
  pages={999--1016},
  year={2022},
  publisher={IEEE}
}

@inproceedings{Bunne,
  title={Learning generative models across incomparable spaces},
  author={Bunne, Charlotte and Alvarez-Melis, David and Krause, Andreas and Jegelka, Stefanie},
  booktitle={International conference on machine learning},
  pages={851--861},
  year={2019},
  organization={PMLR}
}

@inproceedings{Chen,
  title={Semidefinite Relaxations of the {G}romov--{W}asserstein Distance},
  author={Junyu Chen and Binh T Nguyen and Shang Hui Koh and Yong Sheng Soh},
  booktitle={Advances in Neural Information Processing Systems 37},
  year={2024}
}

@article{Memoliultra,
  title={The Ultrametric {G}romov--{W}asserstein Distance},
  author={M{\'e}moli, Facundo and Munk, Axel and Wan, Zhengchao and Weitkamp, Christoph},
  journal={Discrete \& Computational Geometry},
  volume={70},
  number={4},
  pages={1378--1450},
  year={2023},
  publisher={Springer}
}

@Inbook{Vavasis2009,
author="Vavasis, Stephen A.",
editor="Floudas, Christodoulos A.
and Pardalos, Panos M.",
title="Complexity Theory: Quadratic Programming",
bookTitle="Encyclopedia of Optimization",
year="2009",
publisher="Springer US",
address="Boston, MA",
pages="451--454",
isbn="978-0-387-74759-0",
doi="10.1007/978-0-387-74759-0_81",
url="https://doi.org/10.1007/978-0-387-74759-0_81"
}

@article{Chowdhury,
  title={The {G}romov--{W}asserstein distance between networks and stable network invariants},
  author={Chowdhury, Samir and M{\'e}moli, Facundo},
  journal={Information and Inference: A Journal of the IMA},
  volume={8},
  number={4},
  pages={757--787},
  year={2019},
  publisher={Oxford University Press}
}

@inproceedings{Peyre,
  title={Gromov--{W}asserstein averaging of kernel and distance matrices},
  author={Peyr{\'e}, Gabriel and Cuturi, Marco and Solomon, Justin},
  booktitle={International conference on machine learning},
  pages={2664--2672},
  year={2016},
  organization={PMLR}
}

@article{Demetci,
  title={{SCOT}: single-cell multi-omics alignment with optimal transport},
  author={Demetci, Pinar and Santorella, Rebecca and Sandstede, Bj{\"o}rn and Noble, William Stafford and Singh, Ritambhara},
  journal={Journal of Computational Biology},
  volume={29},
  number={1},
  pages={3--18},
  year={2022},
  publisher={Mary Ann Liebert, Inc., publishers 140 Huguenot Street, 3rd Floor New~…}
}

@article{Kravtsova,
  title={Scalable {G}romov--{W}asserstein Based Comparison of Biological Time Series},
  author={Kravtsova, Natalia and McGee II, Reginald L. and Dawes, Adriana T.},
  journal={Bulletin of Mathematical Biology},
  volume={85},
  number={8},
  pages={77},
  year={2023},
  publisher={Springer}
}

@inproceedings{Scetbon,
  title={Linear-time {G}romov--{W}asserstein distances using low rank couplings and costs},
  author={Scetbon, Meyer and Peyr{\'e}, Gabriel and Cuturi, Marco},
  booktitle={International Conference on Machine Learning},
  pages={19347--19365},
  year={2022},
  organization={PMLR}
}

@article{Riahi,
  title={Empot: Partial alignment of density maps and rigid body fitting using unbalanced gromov-wasserstein divergence},
  author={Riahi, Aryan Tajmir and Zhang, Chenwei and Chen, James and Condon, Anne and Duc, Khanh Dao},
  journal={arXiv preprint arXiv:2311.00850},
  year={2023}
}

@article{kravtsova2025,
  title={k-sample inference via multimarginal optimal transport},
  author={Kravtsova, Natalia},
  journal={Electronic Journal of Statistics},
  volume={19},
  number={2},
  pages={4356--4400},
  year={2025},
  publisher={The Institute of Mathematical Statistics and the Bernoulli Society}
}

@misc{kravtsova2026note,
      title={Note on computational complexity of the Gromov-Wasserstein distance}, 
      author={Natalia Kravtsova},
      year={2026},
      eprint={2408.06525},
      archivePrefix={arXiv},
      primaryClass={stat.ML},
      url={https://arxiv.org/abs/2408.06525},
}

@inproceedings{ryner2023,
  author    = {Ryner, Martin and Kronqvist, Jan and Karlsson, Johan},
  title     = {Globally solving the Gromov-Wasserstein problem for point clouds in low dimensional Euclidean spaces},
  booktitle = {Advances in Neural Information Processing Systems},
  volume    = {36},
  pages     = {7930--7946},
  year      = {2023},
  publisher = {Neural Information Processing Systems Foundation, Inc. (NeurIPS)}
}

@article{wu2026,
  title={Multimodal Alignment Through Joint Kernel Entropic Gromov--Wasserstein Optimal Transport},
  author={Wu, Yixuan Florence and Zhu, Yilun and Shi, Naichen},
  journal={arXiv preprint arXiv:2608.04234},
  year={2026}
}

@article{han2025,
  title={Covariance alignment: from maximum likelihood estimation to Gromov--Wasserstein},
  author={Han, Yanjun and Rigollet, Philippe and Stepaniants, George},
  journal={SIAM Journal on Mathematics of Data Science},
  volume={7},
  number={3},
  pages={1491--1513},
  year={2025},
  publisher={SIAM}
}

@article{hohmeier2026,
  title={$ k $-Nearest Neighbors in Gromov--Wasserstein Space},
  author={Hohmeier, Kaitlyn and Fraiman, Nicolas and Moosmueller, Caroline},
  journal={arXiv preprint arXiv:2606.10295},
  year={2026}
}

@article{weitkamp2024,
  title={Distribution of Distances based Object Matching: Asymptotic Inference},
  author={Weitkamp, Christoph Alexander and Proksch, Katharina and Tameling, Carla and Munk, Axel},
  journal={Journal of the American Statistical Association},
  volume={119},
  number={545},
  pages={538--551},
  year={2024},
  publisher={Taylor \& Francis},
  doi={10.1080/01621459.2022.2127360}
}

@article{otcrossdomain2025,
  title={Optimal Transport based Cross-Domain Integration for Heterogeneous Data},
  author={Tang, Shulei and Zhang, Zheng and Zhou, Huijuan and Dunson, David B. and Ma, Rong},
  journal={Journal of the American Statistical Association},
  volume={120},
  number={551},
  pages={542--555},
  year={2025},
  publisher={Taylor \& Francis},
  doi={10.1080/01621459.2025.2540653}
}

@article{zhang2024,
  title={Gromov--Wasserstein distances: Entropic regularization, duality and sample complexity},
  author={Zhang, Zhengxin and Goldfeld, Ziv and Mroueh, Youssef and Sriperumbudur, Bharath K.},
  journal={The Annals of Statistics},
  volume={52},
  number={4},
  pages={1677--1707},
  year={2024},
  publisher={Institute of Mathematical Statistics},
  doi={10.1214/24-AOS2406}
}

@phdthesis{vayer2020,
  author  = {Vayer, Titouan},
  title   = {A contribution to Optimal Transport on incomparable spaces},
  school  = {Universit\'e Bretagne Sud},
  year    = {2020},
  eprint  = {2011.04447},
  archivePrefix = {arXiv},
  primaryClass  = {stat.ML},
  note    = {arXiv:2011.04447}
}

@article{steinkamp2026,
  title={Gromov-Wasserstein Barycenter Surrogates: Statistical Methodology, Distributional Limits and Applications},
  author={Steinkamp, Florian and Rodr{\'\i}guez, Luis-Alberto and Otte, Jan Victor and Munk, Axel},
  journal={arXiv preprint arXiv:2609.11670},
  year={2026}
}

@article{oliver2026,
  title={A Unified Geometric Framework for Developmental Analysis of Spatial Transcriptomic Data},
  author={Oliver, Mary Chriselda Antony and Hohmeier, Kaitlyn and Tran, Tuyen and Castillo, Alejandra and Moosm{\~A}{\v{z}}ller, Caroline and Li, Shiying},
  journal={arXiv preprint arXiv:2608.15306},
  year={2026}
}
